\documentclass[11pt, english]{amsart}
\usepackage{amssymb,amsthm,amsmath,amstext, color}
\usepackage{mathrsfs}
\usepackage{bm}
\usepackage[utf8]{inputenc}
\usepackage{xcolor}
\usepackage
{todonotes}
\usepackage{colortbl}

\usepackage{mathrsfs,upgreek}
\usepackage{appendix}
\usepackage{hyperref}

\input xy
\xyoption{all}

\usepackage[only,mapsfrom,heavycircles,lightning]{stmaryrd}
\usepackage{calc}
\usepackage{paralist}

\makeatletter
\ifnum\@ptsize=0  \fi \ifnum\@ptsize=2  \fi 

\makeatletter
\newcommand{\arkfamily}{\fontencoding{U}\fontfamily{ark}\selectfont}
\newcommand{\ark@sym}[1]{{\arkfamily\symbol{#1}}}
\newcommand{\leftthumbsup}{\ark@sym{'125}}
\newcommand{\smallpencil}{\ark@sym{'120}}

\newcounter{genus}
\newcounter{gagid}

\@addtoreset{gagid}{section}

\makeatother

\theoremstyle{plain}
\newtheorem{theorem}{Theorem}
\newtheorem*{conjecture*}{Conjecture}
\newtheorem{prop}[theorem]{Proposition}
\newtheorem{lemma}[theorem]{Lemma}
\newtheorem{coro}[theorem]{Corollary}
\newtheorem{remark}[theorem]{Remark}

\def\Hom{{\mathrm{Hom}}}
\def\LG{{\mathbf{LG}}}
\def\G{{\mathbf{G}}}
\def\CC{{\mathbb{C}}}
\def\RR{{\mathbb{R}}}
\def\AA{{\mathbb{A}}}
\def\OO{{\mathbb{O}}}
\def\HH{{\mathbb{H}}}
\def\PP{{\mathbb{P}}}
\def\QQ{{\mathbb{Q}}}\def\ZZ{{\mathbb{Z}}}
\def\SS{{\mathbb{S}}}
\def\cO{{\mathcal{O}}}

\def\cM{{\mathcal{M}}}

\def\ra{{\rightarrow}}
\def\lra{{\longrightarrow}}

\def\fc{{\mathfrak c}}\def\fe{{\mathfrak e}}\def\ff{{\mathfrak f}}
\def\fg{{\mathfrak g}}
\def\fso{\mathfrak{so}}
\def\fsp{\mathfrak{sp}}\def\fh{{\mathfrak h}}
\def\fsl{\mathfrak{sl}}

\begin{document}

\title{$F_4$ and desmic quartic surfaces}

\author{Laurent Manivel}
\address{Toulouse Mathematics Institute \\ CNRS/Toulouse University \\  France}
\email{manivel@math.cnrs.fr}
\date{}

\maketitle

\begin{abstract}
    The desmic pencil of quartic surfaces is part of a beautiful, but mostly forgotten chapter of the classical theory of algebraic surfaces: it is the only non-degenerate pencil of surfaces in $\PP^3$ containing at least three completely reducible members. We observe in this note that it is closely related to the Weyl group of the root system $F_4$, and can be recovered from a series of symmetric spaces deduced from the exceptional Lie algebras. We discuss the main properties of the pencil from this Lie theoretic point of view.  
\end{abstract}

\section{Introduction} 

The desmic pencil was named this way by Stephanos in 1879 \cite{stephanos}: it is the unique pencil of quartic surfaces, up to projectivity, containing three tetrahedra with no common face. It was studied in depth by Humbert \cite{humbert} and is devoted a chapter of Jessop's classical treatise on quartic surfaces \cite{jessop}. The general desmic surface has $12$ singular points (the desmic points) and contains $16$ lines and $16$ conics. Its symmetry group is the Weyl group of the root system $D_4$, but we observe that the symmetry group of the full desmic pencil is bigger: it is the Weyl group of $F_4$. 

In fact we will proceed in the reverse direction. In a recent series of works \cite{genre2, genre3, yingqi, gopel}, we studied gradings of simple complex Lie algebras as a powerful tool to construct interesting geometric objects, following the impulse given by \cite{bhargava-ho} and \cite{gsw}. The main tools for dealing with such gradings were devised by Vinberg \cite{vinberg}, whose tables show that the Weyl group of $F_4$ is the so-called little Weyl group of a nice series of $\ZZ_2$-gradings of the exceptional Lie algebras. Typically, there is $\ZZ_2$-grading of $\fe_6$ with one piece 
isomorphic to $\CC^2\otimes\wedge^3\CC^6$, hence connected to lines in $\PP(\wedge^3\CC^6)$. 
This projective space contains an invariant quartic hypersurface \cite[Chapter 3]{donagi}, and such a line therefore defines an elliptic curve: this process, considered in \cite[6.4.2]{bhargava-ho}, was first suggested by Segre \cite{segre}. Applying Vinberg theory, we construct in $\CC^2\otimes\wedge^3\CC^6$ a four-dimensional Cartan subspace (a generalization of the Cartan subalgebra). Restricting the invariant quartic, we get a $W(F_4)$-invariant pencil of quartic surfaces in $\PP^3$, and this turns out to be the desmic pencil.
This is another confirmation of the author's fantasy that part of the most beautiful objects of classical algebraic geometry are, in some sense, derived from the exceptional Lie algebras: including the Segre cubic,  Kummer surfaces, Coble hypersurfaces, and even cubic fourfolds and Debarre-Voisin fourfolds \cite{im-tits}.

Interestingly, the three tetrahedra in the desmic pencil are in correspondence with the twelve long positive roots of $F_4$, that one can uniquely split into three groups of mutually orthogonal roots. One can do the same with the short roots, which yields another copy of the desmic pencil, classically known as conjugate to the first one (its vertices are the $12$ desmic points).
From this point of view, the two desmic pencils are instances of (projectivized) Macdonald representations \cite{macdonald}, constructed from subroot systems of $F_4$ of type $A_1^4$. It is worth noticing that little is known about linear systems of hypersurfaces containing many completely reducible members. For pencils, some progress have been made \cite{py, yuz}; there is a combinatorial version called the desmic conjecture (see e.g. \cite{borwein});  but the general case seems wide open. Macdonald representations provide nice examples of such linear systems. In \cite{liu-manivel}, we found a $\PP^4$ (resp. a $\PP^8$ and a $\PP^4$) of hypersurfaces in $\PP^3$ containing $15$ tetrahedra (resp. $45$ octahedra and $15$ dodecahedra); one can wonder whether these examples are extremal and isolated. 

Closely related to the desmic pencil is the Hesse pencil of plane cubics, which can also be derived from Vinberg's theory. As well as the Hesse pencil parametrizes elliptic curves, it was shown in 
\cite{humbert} that the quartic surfaces of the desmic pencil are biratonal to Kummer surfaces of squares of elliptic curves. We connect this wonderful observation to 
our approach through Vinberg's theory by checking that the elliptic curve in question is, as one could expect, the same as  the one considered by Segre and Bhargava-Ho; this is the only original result of this modest note. 

\medskip \noindent 

{\it Acknowledgements.} We warmly thank Willem de Graaf for his very useful hints, and Igor Dolgachev for his kind comments. Support is acknowledged from the ANR project FanoHK ANR-20-CE40-0023.

\section{$F_4$ and its root system}

A convenient way to describe the root system of $F_4$ is to use the fact that $\fso_9$ embeds in $\ff_4$ as a simple subalgebra of maximal rank. More precisely, there is a 
$\ZZ_2$-grading of $\ff_4$ given by
$$\ff_4 = \fso_9\oplus\Delta,$$
where $\Delta$ is the sixteen dimensional spin representation of $\fso_9$. As a consequence, the standard Cartan subalgebra of $\fso_9$ is also a Cartan subalgebra of $\ff_4$, with respect to which the roots of the latter are 
$$\pm \epsilon_i\pm \epsilon_j, \quad \pm \epsilon_k, \quad 
\frac{1}{2}(\pm \epsilon_1\pm \epsilon_2\pm \epsilon_3\pm \epsilon_4),$$
where $\epsilon_1, \epsilon_2, \epsilon_3, \epsilon_4$ is the usual notation for an
orthogonal basis of the dual of the Cartan algebra. 
The first group gives the $24$ long roots, the second and third group give the $8+16$ short roots
\cite{bourbaki}. In order to select a set positive roots, we can fix a generic real linear form, say $\ell = \ell_1\epsilon_1^\vee+\ell_2\epsilon_2^\vee+\ell_3\epsilon_3^\vee+\ell_4
\epsilon_4^\vee$; the positive roots are those on which $\ell$ takes positive values, and we deduce the simple roots as the indecomposable positive roots; fixing for example 
$\ell_1\gg\ell_2\gg\ell_3\gg\ell_4\gg1$, we obtain the simple roots 
\begin{equation}\label{simple}
\alpha_1=\epsilon_2-\epsilon_3, \quad \alpha_2=\epsilon_3-\epsilon_4, \quad 
\alpha_3=\epsilon_4, \quad \alpha_4=\frac{1}{2}(\epsilon_1-\epsilon_2-\epsilon_3-\epsilon_4).
\end{equation}

\begin{lemma}
Any long (resp. short) root is orthogonal to three long (resp. short) positive roots and six short (resp. long) positive roots.
\end{lemma}

\proof Since the Weyl group acts transitively on long (resp. short) roots, it suffices to check the claim for one long root and one short root. 

The long root $\epsilon_1-\epsilon_2$ is orthogonal to the three long positive roots
 $\epsilon_1+\epsilon_2$ and  $\epsilon_3\pm\epsilon_4$,  and the six short positive roots $\epsilon_3$,   $\epsilon_4$ and  $\frac{1}{2}(\epsilon_1+\epsilon_2\pm \epsilon_3\pm \epsilon_4)$. 
 
The short root $\epsilon_1$ is orthogonal to the three short positive roots $\epsilon_2$, $\epsilon_3$,  $\epsilon_4$ and the six long positive roots $\epsilon_2\pm \epsilon_3$, $\epsilon_2\pm \epsilon_4$, $\epsilon_3\pm \epsilon_4$. \qed 

\medskip An immediate consequence is the following 

\begin{prop}\label{3groups}
There is a unique way to split the twelve long (resp. short) positive roots of $F_4$ into three mutually orthogonal groups of four roots. 
\end{prop}

Explicitly, the three groups of long roots are
$$\begin{array}{llll}
\epsilon_1-\epsilon_2, & \epsilon_1+\epsilon_2, &\epsilon_3-\epsilon_4,& \epsilon_3+\epsilon_4, \\ 
\epsilon_1-\epsilon_3, &\epsilon_1+\epsilon_3, &\epsilon_2-\epsilon_4,& \epsilon_2+\epsilon_4, \\ 
\epsilon_1-\epsilon_4, &\epsilon_1+\epsilon_4, &\epsilon_2-\epsilon_3,& \epsilon_2+\epsilon_3.
\end{array}$$
The three groups of short roots are $
\epsilon_1, \epsilon_2, \epsilon_3, \epsilon_4$ and
$$\begin{array}{ll}
\frac{1}{2}(\epsilon_1+\epsilon_2+\epsilon_3+\epsilon_4), &
\frac{1}{2}(\epsilon_1-\epsilon_2-\epsilon_3+\epsilon_4), \\
\frac{1}{2}(\epsilon_1-\epsilon_2+\epsilon_3-\epsilon_4), &
\frac{1}{2}(\epsilon_1+\epsilon_2-\epsilon_3-\epsilon_4), \\
 \\
\frac{1}{2}(\epsilon_1-\epsilon_2-\epsilon_3-\epsilon_4), &
\frac{1}{2}(\epsilon_1-\epsilon_2+\epsilon_3+\epsilon_4), \\
\frac{1}{2}(\epsilon_1+\epsilon_2-\epsilon_3+\epsilon_4), &
\frac{1}{2}(\epsilon_1+\epsilon_2+\epsilon_3-\epsilon_4).
\end{array}$$

Subsets of mutually orthogonal roots as in Proposition \ref{3groups} are in correspondence with root subsystems of type $A_1^4$. Such root subsystems have been extensively studied and classified \cite{oshima}.

\begin{remark}\label{trans}
The orthogonal matrix 
$$\frac{1}{\sqrt{2}}\begin{pmatrix}
1&-1&0&0 \\ 1&1&0&0 \\ 0&0&1&-1 \\ 0&0&1&1
\end{pmatrix}$$
sends the short roots to multiples of the long roots and conversely. This explains the symmetry between long and short roots in the previous statement.
\end{remark}


\section{$W(F_4)$ as a little Weyl group}

Vinberg showed how to associate to a simple graded Lie algebra a {\it little Weyl group}, generalizing the construction of the usual Weyl group \cite{vinberg}. It is a remarkable fact that, as can be observed in \cite[Table pp. 491--492]{vinberg},  the Weyl group of $F_4$ is the little Weyl group of a series of $\ZZ_2$-gradings, namely  
$$\ff_4 = \fsl_2\times\fsp_6\oplus \CC^2\otimes \wedge^{\langle 3\rangle}\CC^6, $$
$$\fe_6 = \fsl_2\times\fsl_6\oplus \CC^2\otimes \wedge^3\CC^6, $$
$$\fe_7 = \fsl_2\times\fso_{12}\oplus \CC^2\otimes \Delta_+, $$
$$\fe_8 = \fsl_2\times\fe_7\oplus \CC^2\otimes V_{56}. $$

\medskip
Here the fundamental representations $\wedge^{\langle 3\rangle}\CC^6, \wedge^3\CC^6,  \Delta_+, V_{56}$ of $\fsp_6, \fsl_6, \fso_{12}, \fe_7$ have dimensions $6a+8$, for $a=1, 2, 4, 8$. 
This corresponds to the so called sub-exceptional, or  sub-adjoint series, that has been studied by several authors, either from the representation-theoretic \cite{deligne-gross, lm-deligne} or the geometric point of view \cite{freudenthal,lm-freudenthal,buczynski,yingqi}. 

\smallskip
As in the classical theory, the little Weyl group comes with an action on the Cartan subspace (generalizing the usual action of the Weyl group on the Cartan algebra); for each grading of the sub-adjoint series, the Cartan subspace as dimension four and the action of the little Weyl group is equivalent to the action of $W(F_4)$  on a Cartan subalgebra of $\ff_4$. In particular the associated symmetric spaces all have rank four, and the corresponding algebras of invariants are all freely generated by homogeneous polynomials of degrees $2,6, 8, 12$. 

Cartan spaces can be computed explicitly. One can check that the sub-adjoint series of  Lie algebras 
$$\fsp_6\subset \fsl_6\subset \fso_{12}\subset\fe_7$$
comes with natural embeddings of representations 
$$\wedge^{\langle 3\rangle}\CC^6\subset \wedge^3\CC^6\subset  \Delta_+\subset V_{56},$$
and that one can choose a Cartan subspace in $ \CC^2\otimes \wedge^{\langle 3\rangle}\CC^6,$ that can serve as a Cartan subspace for the whole series.

Explicitly, choose a basis $f_1,f_2$ of $\CC^2$ and a basis $e_1,\ldots ,e_6$ of $\CC^6$; denote $e_{ijk}=e_i\wedge e_j\wedge e_k$. We claim that a Cartan subspace $\fc$ can be chosen to be generated by
$$\begin{array}{rcl}
c_1 & = f_1\otimes e_{135} + f_2\otimes e_{246}, \\  
c_2 & = f_1\otimes e_{146} + f_2\otimes e_{235}, \\  
c_3 & = f_1\otimes e_{245} + f_2\otimes e_{136}, \\  
c_4 & = f_1\otimes e_{236} + f_2\otimes e_{145}.
\end{array}$$
Here $e_{135}$ and the other trivectors can be seen as vectors in $\wedge^3\CC^6$, or also 
in $ \wedge^{\langle 3\rangle}\CC^6$, which simply mean that that they correspond to three-spaces that are isotropic with respect to the same symplectic form  $$\omega =e_1^\vee\wedge e_2^\vee+e_3^\vee\wedge e_4^\vee+e_5^\vee\wedge e_6^\vee.$$
In the sequel we represent elements of $W(F_4)$ as matrices encoding their action on $\fc$ in the basis above. 

\begin{prop}\label{generators}
The little Weyl group $W(F_4)$ is generated by the following transformations:
\end{prop}


$$\begin{array}{ll}
s_1 = \begin{pmatrix} 1&   0&   0&  0 \\ 0& 0&  1&    0 \\
 0&   1& 0&  0 \\  0&   0&   0&   1 \end{pmatrix}, & \qquad 
s_2 = \begin{pmatrix}    1& 0&  0&  0 \\ 0& 1&   0&   0 \\
  0&   0&  0& 1 \\  0&   0& 1&  0   \end{pmatrix}, \\
  & \\
   s_3=\begin{pmatrix}   1& 0&   0&   0 \\   0&  1&   0&   0 \\
    0&   0& 1&  0 \\   0&   0&   0&  -1 \end{pmatrix}, & \qquad
s_4= \frac{1}{2}\begin{pmatrix}   1&  1&  1&  1 \\
  1&  1&  -1&  -1 \\
    1&  -1&   1&  -1 \\
   1&  -1&  -1&   1 
   \end{pmatrix}.
   \end{array}$$

\medskip 

\proof First observe that  $s_1,s_2,s_3,s_4$ are hyperplane reflections with respect to 
$\epsilon_2-\epsilon_3$, $\epsilon_3-\epsilon_4$, $\epsilon_4$ and $\frac{1}{2}(\epsilon_1
-\epsilon_2-\epsilon_3-\epsilon_4)$, which is our system (\ref{simple}) of simple roots for $F_4$. As a consequence, they must generate a copy of $W(F_4)$. Since we already know that the little Weyl group is
isomorphic to $W(F_4)$, we just need to check that it contains  $s_1,s_2,s_3,s_4$.

By definition, the little Weyl group is the quotient of the normalizer of $\fc$ in $SL_2\times SL_6$, by its centralizer. Therefore,  we only  need to check that we can lift   $s_1,s_2,s_3,s_4$ to $SL_2\times SL_6$. We provide explicit lifts $\sigma_1, \sigma_2, \sigma_3, \sigma_4$:

$$\begin{array}{l}
\sigma_1 = \begin{pmatrix} -1&0\\ 0&-1\end{pmatrix}\times 
 \begin{pmatrix} 0&0&0&0&1&0 \\ 0&0&0&0&0&1 \\ 0&0&1&0&0&0 \\ 0&0&0&1&0&0 \\
 1&0&0&0&0&0 \\ 0&1&0&0&0&0 \end{pmatrix}, \\
 \\
 \sigma_2 = \begin{pmatrix} -1&0\\ 0&-1\end{pmatrix}\times 
 \begin{pmatrix} 1&0&0&0&0&0 \\ 0&1&0&0&0&0 \\ 0&0&0&0&1&0 \\ 0&0&0&0&0&1 \\ 0&0&1&0&0&0 \\ 0&0&0&1&0&0 
 \end{pmatrix}, \\
 \\
  \sigma_3 = \begin{pmatrix} \zeta &0\\ 0&\zeta^{-1}\end{pmatrix}\times 
 \begin{pmatrix} 1&0&0&0&0&0 \\ 0&-\zeta^2&0&0&0&0 \\ 0&0&-\zeta^2&0&0&0 \\ 0&0&0&1&0&0 \\ 0&0&0&0&\zeta^3&0 \\ 0&0&0&0&0&\zeta 
 \end{pmatrix}, \quad \zeta^4=-1 
  \end{array}$$

 \vfill\pagebreak 
 
$$\begin{array}{l}
 \sigma_4 = \frac{1}{\sqrt{2}}\begin{pmatrix} 1&-1\\ 1&1\end{pmatrix}\times  \frac{1}{\sqrt{2}}
 \begin{pmatrix} 1&-1&0&0&0&0 \\ 1&1&0&0&0&0 \\ 0&0&1&-1&0&0 \\ 0&0&1&1&0&0 \\ 0&0&0&0&1&-1 \\ 0&0&0&0&1&1 
 \end{pmatrix}.
  \end{array}$$
  
  \smallskip
The verification that these are lifts as claimed is left to the reader. \qed 

\medskip Considering the image in $SL_2$ we deduce: 

\begin{coro}\label{binary}
The image in $SL_2$  of the normalizer of the Cartan subspace is a copy if the binary octahedral group.  
\end{coro}

\begin{remark}
The root system of type $F_4$ is sometimes described in terms of quaternions. Indeed, let us formally replace $\epsilon_1, \epsilon_2, \epsilon_3, \epsilon_4$ by the quaternions $1, i, j, k$. Consider the cube whose vertices are   $\pm 1, \pm i, \pm j, \pm k$. The group of unit quaternions that preserve this cube by conjugation is a copy of the binary octahedral group; it consists in the $24$ short roots, and the $24$ long roots, suitably normalized. 
\end{remark}

\section{Some quartic hypersurfaces}

The sub-adjoint series has the form $\fg=\fsl_2\times\fh\oplus\CC^2\otimes W$ where the representations $W$ share the following properties \cite{lm-freudenthal, buczynski}. 
\begin{enumerate}
    \item They are symplectic representations, being endowed with a non-degenerate 
    invariant skew-symmetric form $\Omega$.
    \item They contain a unique invariant hypersurface, of degree four.
    \item The orbit closures of $H=Aut(\fh)$ in $\PP(W)$ are the invariant quartic $Q$, its singular locus $Y$, and the singular locus $X$ of $Y$, which is $H$-homogeneous.
    \item The quartic $Q$ is the tangent variety of $X$ (the union of its tangent spaces), and also the projective dual variety $X^*\subset \PP(W^\vee)\simeq \PP(W)$ (the union of its tangent hyperplanes).
    \item Any point of $\PP(W)$, not belonging to the quartic $Q$, belongs to a unique  bisecant to $X$.
\end{enumerate}

The quartic can be constructed from the symplectic form. 
J.L. Clerc obtained a nice uniform expression by downgrading the  $\ZZ_2$-grading of $\fg$ into a five-step $\ZZ$-grading with $\fg_2\simeq \fg_{-2}\simeq\CC$ and $\fg_1\simeq W$ \cite{clerc}. A more direct statement is the following. 
Recall that $\fh$ is self-dual through the Cartan-Killing quadratic form $K$. 


\begin{prop}\label{quartic}
\begin{enumerate}
\item[]
    \item There is an invariant quadratic  map $P: W\lra\fh$ defined by 
    $$K(P(w),h)=\Omega(hw,w), \quad \forall h\in\fh, w\in W.$$
  \item The invariant hypersurface $Q$ is defined by the quartic equation
  $$K(P(w),P(w))=\Omega(P(w)w,w)=0, \quad w\in W.$$
\end{enumerate}
\end{prop}

\proof We check that there is a unique invariant linear map $S^2W\lra\fh$, and a unique invariant $S^4W\lra\CC$, up to scalar. So any invariant expression, if nonzero, must be the right one. The claim easily follows. \qed

\begin{remark}\label{elliptic}
A generic line $\delta\subset\PP(W)$ cuts the quartic hypesurface $Q$ at four points $q_1,q_2, q_3,q_4$, and the double cover of $\delta$ branched at these four points is a genus one curve $E_\delta$. 

This is one of the constructions of elliptic curves from special representations 
considered in \cite{bhargava-ho} and \cite{thorne}, which ultimately led to major results on the statistics of elliptic curves over $\QQ$. Modular interpretations are studied further in \cite{yingqi}.
\end{remark}

Now consider a generic element $w$ of the Cartan subspace $\fc$. For simplicity we consider $\fc$ as a subspace of $\CC^2\otimes \wedge^3\CC^6$, but the other cases would lead to the same results. 
So we see $w$ as a pencil of trivectors of the form
$$ s(x_1e_{135}+x_2 e_{146}+x_3 e_{245}+x_4 e_{236}) 
+ t(x_1e_{246}+x_2e_{235}+x_3 e_{136}+x_4e_{145}).$$
Our aim is to evaluate the quartic form $Q$ on such a trivector. 
We first compute $P(w)\in\fsl_6$ using Proposition \ref{quartic}; for $\wedge^3\CC^6$, the invariant symplectic form is $\Omega(w,w')=w\wedge w'\in \wedge^6\CC^6$, that we trivialize by fixing a generator, say $e_{123456}$. If $X=e_i^\vee\otimes e_j$, $i\ne j$, the combinatorics of the Cartan subspace is such that $Xw\wedge w=0$ when $(ij)\ne (12),(34),(45)$. The six nonzero terms are 
$$\begin{array}{rcl}
(e_1^\vee\otimes e_2)w\wedge w & = & 2(s^2x_1x_2-t^2x_3x_4), \\
(e_2^\vee\otimes e_1) w\wedge w & = & 2(s^2x_3x_4-t^2x_1x_2), \\
(e_3^\vee\otimes e_4) w\wedge w & = & 2(s^2x_1x_4-t^2x_2x_3), \\
(e_4^\vee\otimes e_3) w\wedge w & = & 2(s^2x_2x_3-t^2x_1x_4), \\
(e_5^\vee\otimes e_6) w\wedge w & = & 2(s^2x_1x_3-t^2x_2x_4), \\
(e_6^\vee\otimes e_5) w\wedge w & = & 2(s^2x_2x_4-t^2x_1x_3).
\end{array}$$
On the other hand, a traceless diagonal matrix $H$ gives 
$$H w\wedge w = -2st\Big( (h_1+h_3+h_5)x_ 1^2+(h_1+h_4+h_6)x_ 2^2+$$
$$\hspace*{5cm} +(h_2+h_4+h_5)x_ 3^2+(h_2+h_3+h_6)x_4^2 \Big).$$
We deduce that $P(w)$ preserves each of the three planes $\langle e_1,e_2\rangle$, $\langle e_3,e_4\rangle$, $\langle e_5,e_6\rangle$; so it is the sum of its restrictions to these three planes, which are given by the three following matrices in 
$\fsl_2$: 
$$\begin{pmatrix}
st(x_3^2+x_4^2-x_1^2-x_2^2) & 2(s^2x_1x_2-t^2x_3x_4) \\
2(s^2x_3x_4-t^2x_1x_2) & st(x_1^2+x_2^2-x_3^2-x_4^2)
\end{pmatrix}$$
$$\begin{pmatrix}
st(x_2^2+x_3^2-x_1^2-x_4^2) & 2(s^2x_1x_4-t^2x_2x_3) \\
2(s^2x_2x_3-t^2x_1x_4) & st(x_1^2+x_4^2-x_2^2-x_3^2)
\end{pmatrix}$$
$$\begin{pmatrix}
st(x_2^2+x_4^2-x_1^2-x_3^2) & 2(s^2x_1x_3-t^2x_2x_4) \\
2(s^2x_2x_4-t^2x_1x_3) & st(x_1^2+x_3^2-x_2^2-x_4^2).
\end{pmatrix}$$

\smallskip
The contribution of the diagonal part is given by the formula above, with $H$ being the diagonal part of $P(w)$; we get 
$$-2s^2t^2\Big( (x_2^2+x_3^2+x_4^2-3x_1^2)x_1^2+(x_1^2+x_3^2+x_4^2-3x_2^2)x_2^2+$$
$$\hspace{3cm}+(x_1^2+x_2^2+x_4^2-3x_3^2)x_ 3^2+(x_1^2+x_2^2+x_3^2-3x_4^2)x_4^2 \Big).$$

Using the formulas above, we obtain the contribution of the off-diagonal part of $P(w)$ as 
$$8\Big( (s^2x_1x_2-t^2x_3x_4)(s^2x_3x_4-t^2x_1x_2)+(s^2x_1x_4-t^2x_2x_3)(s^2x_2x_3-t^2x_1x_4)$$
$$\hspace*{5cm} +(s^2x_1x_3-t^2x_2x_4)(s^2x_2x_4-t^2x_1x_3)\Big). $$

The sum of these two contributions finally yields:

\begin{prop}\label{desmic}
$$Q(w)=24(s^4+t^4)x_1x_2x_3x_4-6s^2t^2\Big(2\sum_{i<j}x_i^2x_j^2-\sum_kx_k^4\Big).$$
\end{prop}

We thus finally obtain the pencil of quartics of the form 
\begin{equation}\label{desmic1}
8Ax_1x_2x_3x_4-B\Big(2\sum_{i<j}x_i^2x_j^2-\sum_kx_k^4\Big).
\end{equation}
This is classically known as the {\it desmic pencil}.

\section{The desmic pencil} 

The desmic pencil of quartic surfaces is characterized by the fact that it contains three tetrahedra \cite{stephanos}, as follows from the expression given by Proposition \ref{desmic}:
\begin{itemize}
    \item $\mathbf{T_1}$: $st=0$, \\ faces $x_i=0$, \\ vertices $[1,0,0,0], [0,1,0,0], [0,0,1,0], [0,0,0,1]$,
    \item $\mathbf{T_2}$: $s^2-t^2=0$, \\
     faces $\epsilon_1x_1+\epsilon_2x_2+\epsilon_3x_3+\epsilon_4x_4=0$, $ \epsilon_1\epsilon_2\epsilon_3\epsilon_4=1$ \\  vertices $[1,1,1,1], [1,1,-1,-1], [1,-1,1,-1], [1,-1,-1,1]$,
    \item $\mathbf{T_3}$: $s^2+t^2=0$, \\ 
    faces $\epsilon_1x_1+\epsilon_2x_2+\epsilon_3x_3+\epsilon_4x_4=0$, $ \epsilon_1\epsilon_2\epsilon_3\epsilon_4=-1$ \\  vertices $[1,1,1,-1], [1,1,-1,1], [1,-1,1,1], [1,-1,-1,-1]$.
\end{itemize}

\begin{remark}
These three tetrahedra were called "desmiques" by Stephanos (from the Greek $\delta\acute{\epsilon}\sigma\mu\eta$, another word for "bundle" or "sheaf") precisely because they belong to a same pencil of quartic surfaces. 
\end{remark}

We recognize in the three tetrahedra the splitting from Proposition \ref{3groups} of the short roots of $\ff_4$. The long roots can be used as well and yield another classical form of the desmic pencil,  as the pencil of quartic surfaces with equations of the form 
\begin{equation}\label{desmic2}
a(y_1^2-y_2^2)(y_3^2-y_4^2)-b(y_1^2-y_3^2)(y_2^2-y_4^2)+c(y_1^2-y_4^2)(y_2^2-y_3^2)=0, 
\end{equation}
where $a+b+c=0$. This can also be written as 
\begin{equation}\label{desmic3}
(b-c)(y_1^2y_2^2+y_3^2y_4^2)+(c-a)(y_1^2y_3^2+y_2^2y_4^2)+(a-b)(y_1^2y_4^2+y_2^2y_3^2)=0. 
\end{equation}
In fact it was observed already in \cite[section 15]{stephanos} that the twelve desmic points are the vertices of three other tetrahedra forming another desmic system, called conjugate to the first one. The explicit computations of \cite[section 5]{desmic-cubics} show that the desmic pencils (\ref{desmic1})  and (\ref{desmic2}) are conjugate. We conclude:

\begin{prop}\label{two-desmic} 
The two splitting of the long and short roots of $\ff_4$ into three tetrahedra yield two conjugate  desmic pencils of quartic surfaces in $\PP^3$. As a consequence, the Weyl group $W(F_4)$ acts as a group of symmetries of these desmic pencils.
\end{prop}

\begin{remark} 
This was clearly known to Shephard and Todd, being briefly mentioned page 278 of \cite{st};
this comes with a reference to \cite{coxeter}, although desmic surfaces do not appear there. However, Coxeter describes $W(F_4)$, rather denoted $[3,4,3]$, as the symmetry group of the $24$-cell, whose vertices can readily be put in correspondence either with the long or the short roots of $F_4$ (see  "Ces\'aro's construction" in \cite[8.2]{coxeter}). 
\end{remark}

A given desmic surface is of course not invariant under  $W(F_4)$, but (in general) under changes of signs of the coordinates, and also double transpositions: these transformations generate a copy of $W(D_4)$. In particular, as  Kummer surfaces,  
desmic surfaces belong to the family of Heisenberg invariant quartics studied in  \cite{eklund} (but they are never Kummer surfaces). A more classical study of these surfaces is given in \cite{edge}, where it is shown among other things that they contain $15$ tetrahedra yielding the famous Klein configuration. 

\begin{remark} 
Taking into account the relation 
$$(y_1^2-y_2^2)(y_3^2-y_4^2)-(y_1^2-y_3^2)(y_2^2-y_4^2)+(y_1^2-y_4^2)(y_2^2-y_3^2)=0,$$
we see that fixing a relation like (\ref{desmic2}), for some given values of $(a,b,c)$, amounts 
to fixing the cross-ratio of $y_1^2, y_2^2, y_3^2, y_4^2$.
\end{remark}

\begin{remark} 
The relationship with root subsystems of $F_4$ of type $A_1^4$ show that the two desmic pencils are incarnations of the so-called Macdonald representations. The general recipe is the following: given a subsystem of a root system, the product of the corresponding linear forms, and the translates of this product by the Weyl group, span an irreducible representation of the latter \cite{macdonald}. In the case at hand, each subsystem has only three different translates (up to scalar), and they generate a two-dimensional Macdonald representation: after projectivization, a desmic pencil. In fact  the generators $s_3$ and $s_4$ act trivially and the representation factors through the symmetric group $S_3$. Taking the two desmic pencils into account yields the classical extension \cite{atlas}
$$1\lra\mu_2^5\lra W(F_4)\lra S_3\times S_3\lra 1.$$
The stabilizer of $[T_1]$ in $W(F_4)$ clearly acts by signed permutations of the $x_i$'s, and is therefore  a copy of $W(D_4)$. Restricting the previous extension we get 
$$1\lra\mu_2^5\lra W(D_4)\lra S_3\lra 1,$$
which is a manifestation of triality.
\end{remark}

\smallskip
Let us now come back to the form (\ref{desmic1}) and recall a few classical properties of the desmic pencil.
\begin{enumerate}
\item The twelve vertices of the three desmic tetrahedra correspond to completely decomposable elements in the Cartan subspace (for example, $s=t$ and $x_1=x_2=x_3=x_4$ yields $(e_1+e_2)\wedge (e_3+e_4)\wedge (e_5+e_6)$. The desmic property is that any edge of one octahedron meets a pair of opposite edges of the second; the intersection points are the $12$ special points $[0,0,1,\pm 1]$ and correspond to lines in $Y$ (the orbit of partially decomposable forms), typically $[0,0,1,1]$ yields $(te_1+se_2)\wedge (e_{36}+e_{45})\in Y$. These points are called {\it desmic points} in \cite{humbert}. Each desmic quartic of the pencil is singular at the $12$ desmic points.
\item 
The base locus of the desmic pencil  is the union of $16$ lines with equations of the form 
$$ x_1=\epsilon_2x_2+\epsilon_3x_3+\epsilon_4x_4=0, \qquad \epsilon_2, \epsilon_3, \epsilon_4=\pm 1.$$
These $16$ lines meet at the $12$ desmic points: each line contains three of them and each desmic point belongs to four lines: hence a $(16_3, 12_4)$ configuration with a $W(F_4)$-symmetry, which is nothing else than the classical {\it Reye configuration} \cite{dolg, man-config}.  Moreover the tangent plane to each desmic surface in the pencil is constant along any of the $16$ lines; so the residual intersection is a conic, one of the $16$ conics contained in each desmic surface 
\item The product $st(s^4-t^4)$ of the quadratic forms associated to the three tetrahedra  is Klein's invariant of the binary octahedral group, in agreement with Corollary \ref{binary}. 
\end{enumerate}

\section{Desmic and Kummer surfaces}

Let $E$  be the elliptic curve obtained as the quotient of $\CC$ by the lattice $\Lambda$ generated by $\omega_1$ and $\omega_2$. Recall that the Weierstrass  $\wp$-function associated to $\Lambda$ is a doubly periodic meromorphic function with order two poles on $\Lambda$, and regular outside $\Lambda$. The half-periods 
$$e_1=\wp(\frac{\omega_1}{2}), \quad e_2=\wp(\frac{\omega_2}{2}), \quad 
e_3=\wp(\frac{\omega_1+\omega_2}{2})$$
are such that $\wp$ satisfies the differential equation
$$(\wp')^2=4(\wp-e_1)(\wp-e_2)(\wp-e_3).$$
On the other hand, the Weierstrass $\sigma$-function is an odd entire function, quasi-periodic with respect to $\Lambda$, and such that $\sigma'/\sigma$ is a primitive of $-\wp$. Traditionally (see e.g. \cite[section 6.2]{lawden}), one defines three other $\sigma$-functions, $\sigma_1$,  $\sigma_2$,  $\sigma_3$, such that 
the following relations are verified \cite[(6.7.10-11-12)]{lawden}:
$$\wp=\Big(\frac{\sigma_1}{\sigma}\Big) ^2+e_1=\Big(\frac{\sigma_2}{\sigma}\Big) ^2+e_2=\Big(\frac{\sigma_3}{\sigma}\Big) ^2+e_3.$$
From these relations, a straightforward computation (see \cite[section 6]{humbert} or \cite[Chap. 2]{jessop}) shows that, letting 
\begin{equation} \label{param}
y_1=\frac{\sigma_1(u)}{\sigma_1(v)}, \quad 
y_2=\frac{\sigma_2(u)}{\sigma_2(v)}, \quad y_3=\frac{\sigma_3(u)}{\sigma_3(v)},
 \quad y_4=\frac{\sigma(u)}{\sigma(v)},
 \end{equation}
for some $u,v\in\CC$, the following identity holds: 
\begin{equation}\label{half-desmic}
(e_1-e_2)(y_1^2y_2^2+y_3^2y_4^2)+
(e_3-e_1)(y_1^2y_3^2+y_2^2y_4^2)+(e_2-e_3)(y_1^2y_4^2+y_2^2y_3^2)=0. 
\end{equation}

As a consequence, (\ref{param}) provides a rational parametrization of the desmic surface with parameters $(a,b,c)=(e_2,e_3,e_1)$. Since 
the $\sigma$ function is odd, this parametrization  is unchanged when $(u,v)\mapsto 
(-u,-v)$, and one deduces that the desmic surface is birational to the Kummer surface of 
$E\times E$ (see also \cite{dolgachev-kondo} for an alternative approach).
\medskip

But what is this elliptic curve? That's an easy guess in our  context.  Indeed a point in $x= [w]\in\PP(\fc)$, outside the three tetrahedra, belongs to a unique desmic surface $\Delta_x$ in the desmic pencil. On the other hand, it defines a line $\delta(x)$ in $\PP(\wedge^3\CC^6)$, that cuts the invariant quartic at four points; hence defines an elliptic curve $E_{\delta(x)}$, as already noticed in Remark \ref{elliptic}.   

\begin{prop} 
The desmic surface $\Delta_x$ is birational to $\mathrm{Kum}(E_{\delta(x)}\times E_{\delta(x)})$.
\end{prop}

\proof The desmic surface to which $x$ belongs has parameters $A,B$ as in (\ref{desmic1}). 
The four points of $\delta(x)$ correspond to the four points $[s,t]$ in $\PP^1$ such that $B(s^4+t^4)=2As^2t^2$. If $[s,t]$ is one of them, the three others are $[s,-t], [t,s], [t,-s]$; computing their cross-ratio, we find e.g. $-(s^2-t^2)^2/4s^2t^2$.  

On the other hand, applying the transformation of Remark \ref{trans} we find a projectively equivalent desmic surface in the conjugate pencil, with equation 
$$4s^2t^2(y_1^2y_2^2+y_3^2y_4^2)+
(s^2-t^2)^2(y_1^2y_3^2+y_2^2y_4^2)-(s^2+t^2)^2(y_1^2y_4^2+y_2^2y_3^2)=0. $$
Comparing with (\ref{half-desmic}), we see that the half-periods $(e_1,e_2,e_3)$ are such that 
$(e_1-e_2,e_3-e_1,e_2-e_3)$ is proportional to $(4s^2t^2,(s^2-t^2)^2,-(s^2+t^2)^2)$. The elliptic curve $E$ with these half-periods is the plane cubic of equations $y^2=4(x-e_1)(x-e_2)(x-e_3)$, hence a double cover of $\PP^1$ branched at $(e_1,e_2,e_3,\infty)$. Computing the cross-ratio of these four points, we get e.g.  
$$\frac{e_3-e_1}{e_2-e_1} = \frac{(s^2-t^2)^2}{-4s^2t^2}.$$
This proves that $E$ and the curve $E_{\delta(x)}$ are isomorphic. \qed

\begin{remark} There is a similar story with the complex reflection group $G_{26}$ of the Shephard-Todd classification, also known as the Hessian group. Indeed, this is the little Weyl group of a series of $\ZZ_3$-gradings of the exceptional Lie algebras. This series is related to the second row of the Tits-Freudenthal magic square, while the sub-adjoint series is related to its third row \cite{lm-freudenthal}. For this series the rank (the dimension of the Cartan subspace) is equal to three, and the related classical geometric object is the Hesse pencil of cubic curves, which contains four triangles and has many similarities with the desmic pencil. Notwistanding the fact the Hesse pencil parametrizes elliptic curves, while the desmic pencil parametrizes (essentially) squares of such elliptic curves, it was observed long ago that the two pencils are intimately related, see e.g. \cite{desmic-cubics}. A  modern exposition of the geometry of the Hesse pencil can be found in \cite{ad}. 
\end{remark}

\section{Complements}

Starting from the desmic family ${\mathcal X}\subset\PP^1\times\PP^3$ one may consider the double cover $X$ of the latter branched along  ${\mathcal X}$. This fourfold has for canonical bundle  the pull-back of $\cO_{\PP^3}(-2)$, and the projection to $\PP^3$ is an elliptic fibration. Of course $X$ is acted on by $W(F_4)$, but it is  singular over the pre-images of the $16$ desmic points. 
By base change over $\PP^1$, one gets a covering $X\lra X_0$ with $X_0$ a  $W(F_4)$-invariant, also singular Fano fourfold.  This is similar to the construction of \cite[section 7]{ad}; in  
their Theorem 8.6 they obtain a K3 surface which is birational to the Kummer surface of the square of an elliptic curve, just as a desmic quartic. This is another way to see the geometry of the Hesse pencil as a degeneration of that of the desmic pencil.  

\medskip


\bibliography{desmic}

\bibliographystyle{amsalpha}

\end{document}